\documentclass[10pt]{amsart}
\usepackage{amsfonts}
\usepackage{amsmath}
\usepackage{amsthm}
\usepackage{amssymb}
\usepackage{latexsym}
\newtheorem{cor}{Corollary}[section]
\newtheorem{lem}{Lemma}[section]

\newtheorem{prop}{Proposition}[section]

\newtheorem{lemma}{Lemma}[section]
\theoremstyle{definition}
\newtheorem{defn}{Definition}[section]
\theoremstyle{definition}
\newtheorem{thm}{Theorem}

\newtheorem*{rem}{Remark}

\newtheorem{definition}[lemma]{Definition}

\theoremstyle{remark}

\numberwithin{equation}{section}
\def\Definition{\begin{definition}}
\def\Def{\begin{definition}}
\def\endefinition{\end{definition}}
\def\endef{\end{definition}}
\def\endeq{\end{eqnarray}}
\def\eneq{\end{eqnarray}}
\def\endeqn{\end{eqnarray*}}
\def\eneqn{\end{eqnarray*}}
\def\eq{\begin{eqnarray}}
\def\eqn{\begin{eqnarray*}}
\def\beq{\begin{eqnarray}}
\def\beqn{\begin{eqnarray*}}
\def\endeq{\end{eqnarray}}
\def\endeqn{\end{eqnarray*}}
\def\qbox#1{\hbox{\quad#1\quad}}
\def\on{\operatorname}

\newcommand{\nc}{\newcommand}
\newcommand{\rnc}{\renewcommand}
\nc{\cal}{\mathcal}
\nc{\goth}{\mathfrak}
\rnc{\bold}{\mathbf}
\renewcommand{\frak}{\mathfrak}
\renewcommand{\Bbb}{\mathbb}

\nc{\Usb}{U'_q(\fb)}

\nc{\tU}{{\tilde U}_q(\ag)}
\nc{\tL}{\tilde L}
\nc{\U}{{U_q}}
\nc{\set}{{\,;\,}}
\nc\Uf{U_q^-(\ag)}
\nc\Ue{U_q^+(\ag)}
\nc{\seps}{\varepsilon^*}
\nc{\sphi}{\varphi^*}
\nc\wt{{\on{wt}}}
\nc\Wt{{\rm Wt}}

\nc{\Cal}{\cal}
\nc{\Xp}[1]{X^+(#1)}
\nc{\Xm}[1]{X^-(#1)}

\nc{\ch}{\mbox{ch}}
\nc{\Z}{{\bold Z}}
\nc{\J}{{\cal J}}
\nc{\C}{{\bold C}}
\nc{\Q}{{\bold Q}}

\nc{\N}{{\Bbb N}}
\nc{\BZ}{{\Z}}
\nc{\BQ}{{\Bbb Q}}

\nc\lan{\langle}
\nc\ran{\rangle}
\nc\bsl{\backslash}
\nc\mto{\mapsto}
\nc\lra{\leftrightarrow}
\nc\hra{\hookrightarrow}
\nc\sm{\smallmatrix}
\nc\esm{\endsmallmatrix}
\nc\sub{\subset}
\nc\ti{\tilde}
\nc\nl{\newline}
\nc\fra{\frac}
\nc\und{\underline}
\nc\ov{\overline}
\nc\ot{\otimes}
\nc\bbq{\bar{\bq}_l}
\nc\bcc{\thickfracwithdelims[]\thickness0}
\nc\ad{\text{\rm ad}}
\nc\Ad{\text{\rm Ad}}
\nc\Hom{\text{\rm Hom}}
\nc\End{\text{\rm End}}
\nc\Ind{\text{\rm Ind}}
\nc\Res{\text{\rm Res}}
\nc\Ker{\text{\rm Ker}}
\rnc\Im{\text{Im}}
\nc\sgn{\text{\rm sgn}}
\nc\tr{\text{\rm tr}}
\nc\Tr{\text{\rm Tr}}
\nc\supp{\text{\rm supp}}
\nc\card{\text{\rm card}}
\nc\bst{{}^\bigstar\!}
\nc\he{\heartsuit}
\nc\clu{\clubsuit}
\nc\spa{\spadesuit}
\nc\di{\diamond}

\nc\al{\alpha}
\nc\bet{\beta}
\nc\ga{\gamma}
\nc\de{\delta}
\nc\ep{\varepsilon}
\nc\eps{\varepsilon}
\nc\io{\iota}
\nc\om{\omega}
\nc\si{\sigma}
\rnc\th{\theta}
\nc\ka{\kappa}
\nc\lam{\lambda}
\nc\la{\lambda}
\nc\ze{\zeta}

\nc\vp{\varpi}
\nc\vt{\vartheta}
\nc\vr{\varrho}

\nc\Ga{\Gamma}
\nc\De{\Delta}
\nc\Om{\Omega}
\nc\Si{\Sigma}
\nc\Th{\Theta}
\nc\La{\Lambda}

\nc\boa{\bold a}
\nc\bob{\bold b}
\nc\boc{\bold c}
\nc\bod{\bold d}
\nc\boe{\bold e}
\nc\bof{\bold f}
\nc\bog{\bold g}
\nc\boh{\bold h}
\nc\boi{\bold i}
\nc\boj{\bold j}
\nc\bok{\bold k}
\nc\bol{\bold l}
\nc\bom{\bold m}
\nc\bon{\bold n}
\nc\boo{\bold o}
\nc\bop{\bold p}
\nc\boq{\bold q}
\nc\bor{\bold r}
\nc\bos{\bold s}
\nc\bou{\bold u}
\nc\bov{\bold v}
\nc\bow{\bold w}
\nc\boz{\bold z}

\nc\ba{\bold A}
\nc\bb{\bold B}
\nc\bc{\bold C}
\nc\bd{\bold D}
\nc\be{\bold E}
\nc\bg{\bold G}
\nc\bh{\bold H}
\nc\bi{\bold I}
\nc\bj{\bold J}
\nc\bk{\bold K}
\nc\bl{\bold L}
\nc\bm{\bold M}
\nc\bn{\bold N}
\nc\bo{\bold O}
\nc\bp{\bold P}
\nc\bq{\bold Q}
\nc\br{\bold R}
\nc\bs{\bold S}
\nc\bt{\bold T}
\nc\bu{\bold U}
\nc\bv{\bold V}
\nc\bw{\bold W}
\nc\bz{\bold Z}
\nc\bx{\bold X}

\nc\ca{\Cal A}
\nc\cb{\Cal B}
\nc\cc{\Cal C}
\nc\cd{\Cal D}
\nc\ce{\Cal E}
\nc\cf{\Cal F}
\nc\cg{\Cal G}
\rnc\ch{\Cal H}
\nc\ci{\Cal I}
\nc\cj{\Cal J}
\nc\ck{\Cal K}
\nc\cm{\Cal M}
\nc\cn{\Cal N}
\nc\co{\Cal O}
\nc\cp{\Cal P}
\nc\cq{\Cal Q}
\nc\car{\Cal R}
\nc\cs{\Cal S}
\nc\ct{\Cal T}
\nc\cu{\Cal U}
\nc\cv{\Cal V}
\nc\cz{\Cal Z}
\nc\cx{\Cal X}
\nc\cy{\Cal Y}

\nc\e[1]{E_{#1}}
\nc\ei[1]{E_{\delta - \alpha_{#1}}}
\nc\esi[1]{E_{s \delta - \alpha_{#1}}}
\nc\eri[1]{E_{r \delta - \alpha_{#1}}}
\nc\ed[2][]{E_{#1 \delta,#2}}
\nc\ekd[1]{E_{k \delta,#1}}
\nc\emd[1]{E_{m \delta,#1}}
\nc\erd[1]{E_{r \delta,#1}}

\nc\ef[1]{F_{#1}}
\nc\efi[1]{F_{\delta - \alpha_{#1}}}
\nc\efsi[1]{F_{s \delta - \alpha_{#1}}}
\nc\efri[1]{F_{r \delta - \alpha_{#1}}}
\nc\efd[2][]{F_{#1 \delta,#2}}
\nc\efkd[1]{F_{k \delta,#1}}
\nc\efmd[1]{F_{m \delta,#1}}
\nc\efrd[1]{F_{r \delta,#1}}

\nc\fa{\frak a}
\nc\fb{\frak b}
\nc\fc{\frak c}
\nc\fd{\frak d}
\nc\fe{\frak e}
\nc\ff{\frak f}
\nc\g{\frak g}
\nc\fg{\frak g}
\nc\fh{\frak h}
\nc\fj{\frak j}
\nc\fk{\frak k}
\nc\fl{\frak l}
\nc\fm{\frak m}
\nc\fn{\frak n}
\nc\fo{\frak o}
\nc\fp{\frak p}
\nc\fq{\frak q}
\nc\fr{\frak r}
\nc\fs{\frak s}
\nc\ft{\frak t}
\nc\fu{\frak u}
\nc\fv{\frak v}
\nc\fz{\frak z}
\nc\fx{\frak x}
\nc\fy{\frak y}

\nc\fA{\frak A}
\nc\fB{\frak B}
\nc\fC{\frak C}
\nc\fD{\frak D}
\nc\fE{\frak E}
\nc\fF{\frak F}
\nc\fG{\frak G}
\nc\fH{\frak H}
\nc\fJ{\frak J}
\nc\fK{\frak K}
\nc\fL{\frak L}
\nc\fM{\frak M}
\nc\fN{\frak N}
\nc\fO{\frak O}
\nc\fP{\frak P}
\nc\fQ{\frak Q}
\nc\fR{\frak R}
\nc\fS{\frak S}
\nc\fT{\frak T}
\nc\fU{\frak U}
\nc\fV{\frak V}
\nc\fZ{\frak Z}
\nc\fX{\frak X}
\nc\fY{\frak Y}

\nc\tfi{\ti{\Phi}}
\nc\te{\tilde e}
\nc\tf{\tilde f}

\nc\bF{\bold F}
\rnc\bol{\bold 1}

\nc\ua{{_\ca\bold U}}

\nc\qinti[1]{[#1]_i}
\nc\q[1]{[#1]_q}
\nc\xpm[2]{E_{#2 \delta \pm \alpha_#1}}  
\nc\xmp[2]{E_{#2 \delta \mp \alpha_#1}}
\nc\xp[2]{E_{#2 \delta + \alpha_{#1}}}
\nc\xm[2]{E_{#2 \delta - \alpha_{#1}}}
\nc\hik{\ed{k}{i}}
\nc\hjl{\ed{l}{j}}
\nc\qcoeff[3]{\left[ \begin{smallmatrix} {#1}& \\ {#2}& \end{smallmatrix}
\negthickspace \right]_{#3}}
\nc\qi{q}
\nc\qj{q}

\nc\ufdm{{_\ca\bu}_{\rm fd}^{\le 0}}

\nc\isom{\cong} 
\nc{\isoto}{\mathrel{\longrightarrow{\kern-18pt\raise3.5pt\hbox{$\sim$}}}}

\nc{\pone}{{\Bbb C}{\Bbb P}^1}
\nc{\pa}{\partial}

\nc{\F}{{\cal F}}
\nc{\Sym}{{\goth S}}
\nc{\A}{{\cal A}}
\nc{\arr}{\rightarrow}
\nc{\larr}{\longrightarrow}

\nc{\ri}{\rangle}
\nc{\lef}{\langle}
\nc{\W}{{\cal W}}
\nc{\uqatwoatone}{{U_{q,1}}(\su)}
\nc{\uqtwo}{U_q(\goth{sl}_2)}
\nc{\dij}{\delta_{ij}}
\nc{\divei}{E_{\alpha_i}^{(n)}}
\nc{\divfi}{F_{\alpha_i}^{(n)}}
\nc{\Lzero}{\Lambda_0}
\nc{\Lone}{\Lambda_1}
\nc{\ve}{\varepsilon}
\nc{\vphi}{\varphi}
\nc{\phioneminusi}{\Phi^{(1-i,i)}}
\nc{\phioneminusistar}{\Phi^{* (1-i,i)}}
\nc{\phii}{\Phi^{(i,1-i)}}
\nc{\Li}{\Lambda_i}
\nc{\Loneminusi}{\Lambda_{1-i}}
\nc{\vtimesz}{v_\ve \otimes z^m}

\nc{\asltwo}{\widehat{\goth{sl}_2}}
\nc\ag{\widehat{\goth{g}}}  
\nc\teb{\tilde E_\boc}
\nc\tebp{\tilde E_{\boc'}}

\begin{document}

\title{Crystal structure of level zero \break extremal weight modules}

\author[J. ~Beck]{Jonathan Beck} 

\address{Department of Mathematics,
Bar Ilan University, 52900 Ramat Gan, Israel}
\email{beck@macs.biu.ac.il}

\begin{abstract}
   We consider the crystal structure of the level zero extremal weight
 modules $V(\lambda)$ \cite{K1} using the crystal base of the quantum
 affine algebra constructed in \cite{BCP}.  This approach yields an
 explicit form for extremal weight vectors in the $U^-$ part of each
 connected component of the crystal, which are given as Schur
 functions in the imaginary root vectors.  We show the map
 $\Phi_\lambda$ (\cite[\S 13]{K2}) induces a correspondence between
 the global crystal base of $V(\lam)$ and elements $s_{\boc_0}(z^{-1})
 G(b), b \in B_0(\U[z^{\pm 1}_{i,k}]u')$.
\end{abstract}

\maketitle
\section{Introduction.}

This paper arises from the study of the modified quantum algebra
$\tilde U_q(\g) = \oplus_{\lam\in P}\U(\g) a_\la$ and its associated
crystal structure.  For $\g$ a simple Lie algebra the crystal of
$\tilde U_q(\g)$ has a Peter-Weyl type decomposition and is isomorphic
to the crystal of the quantum coordinate ring $\oplus_{\lambda \in
P^+} V(\lambda) \otimes V(-\lambda)$, where $P^+$ is the set of
dominant weights of $\g$, and $V(\lambda)$ denotes the irreducible
highest weight module of weight $\lambda.$

This decomposition fails when $\ag$ is the affine Lie algebra
associated to $\g$. In fact, the crystal of $\tU$ naturally decomposes
into pieces according to the level of $\lambda \in P$, $B(\tU) =
B(\tU)_+ \oplus B(\tU)_0 \oplus B(\tU)_-,$ and $B(\tU)_\pm$ again have
Peter--Weyl type decompositions.  However, as level zero weights are
not Weyl group conjugates of dominant weights, a similar analysis is
impossible for $B(\tU)_0$.

In studying the level zero part of the crystal, Kashiwara introduced
(\cite{K1}) the extremal weight modules $V(\lambda)$.  These are a
natural variation on highest weight modules, since when $\lambda \in
P^+$, $V(\lambda)$ is the usual irreducible highest weight module
generated by $v_\lambda$.  For each $\lambda = \sum_i \lambda_i
\Lambda_i \in P$, $V(\lambda)$ has a basis given by a subset of
``*-extremal'' elements of the global crystal base of $\U(\g) a_\la$.
The term extremal refers to the fact that the usual highest weight
relations, $E_i v_\lambda = 0, \, F_i^{(\lambda_i + 1)} v_\lambda =
0$, are replaced by the more general condition that $E_i v_\lambda
=0,\, F_i^{(\lambda_i + 1)} v_\lambda = 0,\, \text{ if } \lambda_i \ge
0$, and $F_i v_\lambda = 0, \, E_i^{(-\lambda_i + 1)} v_\lambda = 0,
\text{ if } \lambda_i \le 0,$ as well as new relations determined by
the Weyl group.

We study the crystal structure of $V(\lambda)$ for $\lambda = \sum_i
\lambda_i \varpi_i, \, \lambda_i \ge 0$, where $\varpi_i \in P_0^+$
are the fundamental level zero weights.  Kashiwara proves \cite{K2}
that the $V(\varpi_i)$ are affinizations of certain finite dimensional
$U'_q$--modules which have global crystal bases.  In \cite[\S 13]{K2},
Kashiwara conjectures a description of the crystal structure of
$V(\lambda)$ for arbitrary $\lambda \in P^+_0$ in terms of an
$m$--fold tensor product of the crystal of $V(\varpi_i)$ and Schur
functions.  These conjectures also imply a Peter--Weyl type decomposition
of $B(\tU)$, which will appear in forthcoming work. 

The purpose of this note is to verify these conjectures in the
symmetric untwisted case using the crystal base of $B(\tU)$
constructed in \cite{BCP}.  The key property of this basis is that it
too contains Schur functions naturally, and each component of the
crystal of $V(\lambda)$ contains an extremal element corresponding to
one of these Schur functions.

While this paper was in preparation, Nakajima \cite{N} released a
preprint where the same results are obtained.  The primary difference
in the proofs is that here we avoid using an explicit description of
Lusztig's braid group action on extremal weight modules. We obtain the
key property of connected components of the crystal of $V(\lambda)$
using the basis \cite{BCP} directly, and this requires less
calculation.

\vskip12pt\noindent{\bf Acknowledgements.}  I would like to thank
M. Kashiwara for helpful conversations and hospitality at RIMS.

\section{The algebra $U_q(\ag)$ and background.}

In this section, we review very briefly the relevant material on
quantum affine algebras and crystal bases. 

Let $\ag$ be a symmetric untwisted affine Lie algebra over $\Q$ with
Cartan datum $(\tilde I,\cdot)$. Choose a Cartan subalgebra $\ft
\subset \ag$ such that simple roots $\{\alpha_i\}_{i\in \tilde I}$ and
the simple coroots are $\{h_i\}_{i\in \tilde I}$ are linearly
independent and $\dim \ft = |\tilde I| + 1 = \text{rank } \ag + 1$.
Denote by $\langle \,\cdot\,,\,\cdot\,\rangle \colon \ft^*\times \ft
\to \Q$ the canonical pairing.

Fix the root lattice and coroot lattice by
$$Q=\oplus_i\BZ\alpha_i\subset\ft^*\qbox{and} Q^\vee=\oplus_i\BZ
h_i\subset\ft.$$ We assume that $\tilde I=\{0,1,\ldots,n\}$ and that
the index ${0} \in \tilde I$ is such that $(I = \tilde I \setminus
\{{0}\},\cdot)$ is the Cartan datum of the underlying finite type
algebra $\g.$

Set $Q_\pm=\pm\sum_i\BZ_{\ge 0}\alpha_i$ and
$Q_\pm^\vee=\pm\sum_i\BZ_{\ge 0}h_i$.  Let $\delta\in Q_+$ be the unique
element satisfying $\{\lam\in Q\set \lan h_i,\lam\ran=0\hbox{ for
every $i$}\}=\BZ\delta$.  Similarly we define $c\in Q^\vee_+$ by
$\{h\in Q^\vee\set \lan h,\alpha_i\ran=0\hbox{ for every $i$}\}=\BZ
c$.  We write \beq &&\delta=\sum_ia_i\alpha_i\qbox{and}
c=\sum_ia_i^\vee h_i.  \endeq

We choose a weight lattice $P\subset\ft^*$, and $\Lambda_i \in P, i
\in \tilde I$, satisfying \eq &&\text{$\alpha_i\in P$ and $h_i\in P^*$
for any $i\in \tilde I$.}\\ &&\text{For every $i\in \tilde I$, $\lan
h_j,\Lambda_i\ran=\delta_{ji}$.}  \eneq We set $P^0=\{\lam\in
P\set\lan c,\lam\ran=0\}$. For $i = 1, \dots, n$, let $\varpi_i =
\Lambda_i - a_i^\vee \Lambda_0.$ We denote by $P^0_+$ the dominant
level zero weights of the form $\lambda = \sum_{i=1}^n m_i \varpi_i,
m_i \ge 0.$

Denote by $\car^+$ the set of positive roots for $\ag$.  Let
$\car_>$ (resp. $\car_<$) denote the set of $\alpha \in
\car^+$ for which the real part of $\alpha$ is positive (resp.
negative). We identify the set of positive imaginary roots (with
multiplicity) by defining $\car_0 = \{k \delta | k > 0 \} \times \{1,
\dots, n\} = \{ k \delta^{(i)} \ | \ k > 0,\ i = 1, \dots, n\}$.  Then
define the set of {\em positive roots with multiplicity} as
\begin{equation}  \label{rootdecomposition} \car^+ = \car_> \cup \car_0
\cup \car_<.
\end{equation}

\subsection{Quantum affine algebras and crystals}

We denote by $\U = \U(\ag)$ the quantum affine algebra on generators
$F_i$'s, the $E_i$'s and $K_h$ ($h\in P^*$), see \cite{BCP} for
details.  Let us denote by $\Uf$ (resp. $\Ue$) the subalgebra of $\U$
generated by the $F_i$'s (resp. by the $E_i$'s).  Introduce the
divided powers $F_i^{(n)}= F^n_i/{[n]_q !}$, $E_i^{(n)} = E^n_i/{[n]_q
!}$, where $[n]_q = (q^n - q^{-n})/(q -q^{-1})$.  We denote by
$(U_q)_\Q$ the subalgebra of $\U$ generated by the $F_i^{(n)}$'s, the
$E_i^{(n)}$'s ($i\in \tilde I$) and $K_h$ ($h \in P^*$) over
$\Q[q,q^{-1}]$.

 We refer to \cite{K4} for a comprehensive introduction to crystals.
We mention the details relevant to this paper.  A crystal $B$ is
graded by a weight lattice $P_B$ containing simple roots $\alpha_i, i
\in I',$ for some Cartan datum $(I',\cdot)$.  For each $b \in B$ we
denote its weight by $wt(b) \in P_B.$ $B$ also comes with operators
$\tilde e_i, \tilde f_i: B \rightarrow B \sqcup \{0\}, i \in I'$ such
that if $\tilde e_i(b) \neq 0$ (resp. $\tilde f_i(b) \neq 0$),
$wt(\tilde e_i (b)) = wt(b) + \alpha_i$ (resp. $wt(\tilde f_i (b)) =
wt(b) - \alpha_i$).  Given two crystals $B_1, B_2$, their tensor
product, $B_1 \otimes B_2$ denotes the crystal $\{b_1 \otimes b_2 \, |
\, b_1 \in B_1, b_2 \in B_2 \}$ with crystal operators defined in
[ip.cit., \S 7.3].

The algebra $\Uf$ (resp. $\Ue$) has a crystal base (see \cite{K3})
denoted by $(L(\infty),B(\infty))$ (resp.  $(L(-\infty),B(-\infty))$.
The term crystal base refers to the following additional requirements:
$L(\infty)$ is a $\Q[q]$ sublattice of $\Uf$ which is invariant under
crystal operators $\tilde e_i, \tilde f_i, i \in \tilde I$, such that
$B(\infty)$ is a basis of $L(\infty)/qL(\infty)$ and forms a crystal
under the operators induced from $\tilde e_i$ and $\tilde f_i.$ A
similar description holds for $(L(-\infty),B(-\infty))$.  We denote by
$u_\infty$ (resp. $u_{-\infty}$) the unique element of $B(\infty)$
(resp. $B(-\infty)$) of weight $0$.

Let $-$ be the $\Q$-algebra automorphism of $U_q$ sending $q$ to
$q^{-1}$, $K_h$ to $K_{-h}$, and fixing $E_i, F_i.$ Let $M$ be a
$\U$-module with a crystal base $(L,B)$ (see \cite[\S 2]{K3}).
Define a bar involution $-$ on $M$ to be an involution satisfying
$\ov{au}=\ov a \, (\ov u)$ for $a\in\U$ and $u\in M$.  Let $M_\Q$ be a
$(U_q)_\Q$-submodule of $M$ such that 
\eq \ov{M_\Q} = M_\Q, \text{ and } 
u - \ov{u} \in (q-1) M_\Q \text{ for every } u \in M_\Q.  
\eneq Let
$E = L \cap \ov{L} \cap M_\Q$.  If the map $f:E \rightarrow L/qL$ is
an isomorphism of $\Q[q]$ modules then we say $M$ has a global base.
If $b \in B,$ we set $G(b) = f^{-1}(b).$ In this case $\ov{G(b)} =
G(b)$ and $G(b)$ is called the globalization of $b$.

 Let us denote by $\tU$ the modified quantum universal enveloping
algebra $\oplus_{\lam\in P}\U a_\la$ (see \cite{L, K1}).  Then $\tU$
has a crystal base $(\tilde L, B(\tU))$.  Let $T_\lam$ be the crystal
consisting of a single element $t_\lam$ of weight $\wt(t_\lam)=\lam$,
with $\tilde e_i(t_\lam)=\tilde f_i(t_\lam)= 0$.  As a crystal,
$B(\tU)$ is regular and isomorphic to
\begin{equation} \label{tildecrystal}
\bigsqcup_{\lam\in P}B(\infty)\otimes T_\lam\otimes B(-\infty).
\end{equation}

The property of being regular allows one to define a Weyl group $W$
action on the crystal.  For each $i \in \tilde I$, define 
\eqn &&S_{s_i}b=
\begin{cases}
\tf_i^{\lan h_i,\wt(b)\ran}b &\mbox{if $\lan h_i,\wt(b)\ran\ge 0$,}\\
\te_i^{-\lan h_i,\wt(b)\ran}b &\mbox{if $\lan h_i,\wt(b)\ran\le 0$.}
\end{cases}
\endeqn Then by \cite{K1} the $S_{s_i}$ satisfy the defining relations
of the Weyl group.

One of the main results of \cite{BCP} is an explicit construction of a
crystal base for $U_q^\pm(\ag),$ and by \eqref{tildecrystal} one of
$B(\tU).$ We denote by $T_i$ ($=T^{\prime\prime}_{i,1}$ in
\cite[Chapter 37]{L}) the automorphism of $U_q$ corresponding to the
simple reflection $s_i$, $i=0,\ldots ,n$.  For each $w \in W$, the
$T_i$'s define an automorphism $T_w$ of $\U$.  Using these
automorphisms, for each $\alpha \in \car_> \cup \car_<, i = 1 \dots n,
k > 0$, we define root vectors $E_\alpha$, $\tilde P_{i,k} \in
U_q^+(\ag)$ as in \cite{BCP}.

The $\tilde P_{i,k} \, (1 \le i \le n, k>0)$ are used to construct a
basis of the imaginary parts of $(\U)^+_\Q$ as follows. Let $\boc_0$ be
an $n$--tuple of partitions $(\rho^{(1)}, \rho^{(2)}, \dots,
\rho^{(n)})$ where each $\rho^{(i)} = (\rho^{(i)}_1 \ge \rho^{(i)}_2
\ge \dots )$.  For a partition $\rho$, denote by $\rho'$ its transpose.
For each $i$, define Schur functions in the $\tilde
P_{i,k}$ by
\begin{equation*}
S_{\rho^{(i)}} = \det(\tilde{P}_{i, {\rho'}^{(i)}_k-k+m})_{1\le k,m\le t},
\end{equation*}
where $t\ge l({\rho'}^{(i)})$. This puts the $\tilde P_{i,k}$ in the role of
elementary symmetric functions.  Denote the product over $i = 1,
\dots, n$ of $S_{\rho^{(i)}}$ by
\begin{equation*} 
S^+_{\boc_0} = \prod_{i=1}^n S_{\rho^{(i)}}. 
\end{equation*} 

\begin{defn} Let $\boc_+ \in \N^{\car_>}$ (resp. $\boc_- \in \N^{\car_<}$) and
$\boc_0$ as above.  Denote by $\boc = (\boc_+, \boc_0, \boc_-).$ Each
$\boc$ indexes a basis element of $(\U)^+_\Q$,
\begin{equation*}  \label{basis}
 B_{\boc}^+ = (E_{\boc_+}) \cdot S^+_{\boc_0} \cdot (E_{\boc_-}),
\end{equation*}
which when specialized at $q=1$ becomes an element of the Kostant
$\Z$-form of $U_{q=1}^+(\ag)$. When we refer to an element of this type,
we will call $B^+_\boc$ purely imaginary if $\boc_+ = \boc_- = 0$.
\end{defn}

\begin{prop} \cite{BCP}
 $(L(-\infty), \overline{B^+_\boc} u_{-\infty} \mod q L(-\infty))$ forms
 a crystal base of $U_q^+.$
\end{prop}

There are two more involutions of $U_q$ which we refer to: Let $\psi$
be the automorphism of $\U$ which sends $E_i$ to $F_i$, $F_i$ to $E_i$
and $K_h$ to $K_{-h}$.  It gives a bijection $B(\infty)\simeq
B(-\infty)$.  Let $*$ be the anti--automorphism of $U_q$ which fixes
$E_i$ and $F_i$, and sends $K_h$ to $K_{-h}$.  Restricted to $U_q^+,
*$ gives an bijection $B(\infty)\simeq B(\infty)$.  For the
calculations in this paper we use the crystal bases
$(L(-\infty),(\overline{B_\boc^+})^* \mod qL(-\infty))$ of $U_q^+$ and
$(L(\infty), \psi(\overline{B_\boc^+}) \mod qL(\infty))$ of $U_q^-.$
In what follows, we replace the definitions of root vectors in
\cite{BCP} by those obtained by applying the involutions $*$, $-$, and
$\psi$ as described. So $E_\alpha$ actually refers to
$(\overline{E_\alpha})^*$, $F_\alpha$ actually refers to
$\psi(\overline{E_\alpha})$, $P_{i,-k}$ refers to
$\psi(\overline{P_{i,k}}),$ $B_\boc^- = \psi(\overline{B_\boc^+})$,
etc.  The purpose of applying these involutions is to arrange the root
ordering of the crystal bases to aid the calculations.

With the above remarks, these imaginary root vectors satisfy an
important property following from \cite[Proposition 2.2 and
eq. (4.9)]{BCP}:
\begin{prop} \label{greatid}
\begin{align*} 
 \tag{i} \tilde P_{i,-k} &= F_{\delta - \alpha_i}^{(k)}
  F_{\alpha_i}^{(k)} + q x, \\ \tag{ii} \tilde P_{i,k} &=
  E_{\alpha_i}^{(k)} E_{\delta - \alpha_i}^{(k)} + q x,
\end{align*}
where in \text{(i)}, $x$ is a sum of terms $B_\boc^-$ with
coefficients in $\Z[q]$ where for each term $\boc_- \neq 0$. In
\text{(ii)}, $x$ is a sum of terms $B_\boc^+$ with coefficients in
$\Z[q]$ where for each term $\boc_+ \neq 0.$
\end{prop}

\section{Extremal weight modules}

\subsection{Extremal vectors}
Let $M$ be an integrable $\U$-module.  A vector $u\in M$ of weight
$\lambda\in P$ is called {\em extremal} (see \cite{K1}) if we can
find vectors $\{u_w\}_{w\in W}$ satisfying the following properties:
\eq &&\text{$u_w=u$ for $w=e$,}\\ && \hbox{if $\lan h_i,w\lam\ran\ge
0$, then $e_iu_w=0$ and $f_i^{(\lan h_i,w\lam\ran)}u_w=u_{s_iw}$,}\\
&&\hbox{if $\lan h_i,w\lam\ran\le 0$, then $f_iu_w=0$ and $e_i^{(\lan
h_i,w\lam\ran)} u_w =u_{s_iw}$.}  \eneq If such $\{u_w\}$ exists,
it is unique and $u_w$ has weight $w\lam$.  We denote $u_w$ by
$S_wu$.

Similarly, for a vector $b$ of a regular crystal $B$ with weight $\lam$,
we say that $b$
is an extremal vector
if it satisfies the following similar conditions:
we can find vectors
$\{b_w\}_{w\in W}$ such that
\beq
&&b_w=b\quad\hbox{for $w=e$,}\\
&&\hbox{if $\lan h_i,w\lam\ran\ge 0$ then
$\te_ib_w=0$ and $\tf_i^{\lan h_i,w\lam\ran}b_w=b_{s_iw}$,}\\
&&\hbox{if $\lan h_i,w\lam\ran\le 0$ then
$\tf_i b_w=0$ and $\te_i^{(\lan h_i,w\lam\ran)}b_w=b_{s_iw}$.}
\endeq
Then $b_w$ must be $S_wb$.

For $\lam\in P$, we denote by $V(\lam)$ the $\U$-module generated by
$u_\lam$ with the defining relation that $u_\lam$ is an extremal
vector of weight $\lam$ (see \cite{K1} for details).  It is proved in
\cite{K1} that $V(\lam)$ has a global crystal base
$(L(\lam),B(\lam))$.  Moreover, if $M$ is any integrable $\U$ module
with extremal weight vector $u$ of weight $\lam$, there is a unique
$\U$ homomorphism $\Phi_\lambda :V(\lam) \rightarrow M,$ such that
$\Phi_\lambda(u_\lam) = u.$ On an integral $\tU$ module $M$, we use
the regularized crystal operators $\tilde e_i, \tilde f_i: M
\rightarrow M$ as defined in (\cite[\S 6]{K2}).  In this context,
$\Phi_\lambda$ commutes with the crystal operators $\tilde e_i, \tilde
f_i$.

\subsection{Extremal weight modules $V(\lam)$ for $\lam \in P^0_+$}

Denote by $\boc_0(\lambda)$ the set of $\boc_0 =
(\rho^{(1)},\rho^{(2)}, \dots, \rho^{(n)})$ such that for each $i$,
$\ell(\rho^{(i)}) \le \lam_i = \langle h_i, \lambda \rangle$.  The
following is an important corollary to \cite[Theorem 5.1]{K2}.

\begin{prop} \label{iextremal}
(i) For any $\lam \in P_+^0$, any vector in $B(\lambda)$ is connected
to an extremal weight vector of the form $b_1 \otimes t_\lambda
\otimes u_{-\infty},$ where $b_1$ is purely imaginary with respect to
the crystal base. \\ (ii) Furthermore, all such possible $b_1 \in
B(\infty)$ are given by $S^-_{\boc_0} u_\infty \mod qL(\infty)$ where
$\boc_0 \in \boc_0(\lambda).$
\end{prop}

\begin{proof} (i) By \cite[Theorem 5.1]{K2} any vector is connected to 
an extremal weight vector of the form $b_1 \otimes t_\lambda \otimes
u_{-\infty},$ where $wt(b_1) = -k\delta.$ Using the crystal base to
express $b_1,$ we take $B^-_\boc = b_1 \mod qL(\infty),$ for some
$\boc$.  Assume that $B^-_\boc$ isn't purely imaginary.  Since
$wt(B^-_\boc) = -k\delta$, and $B^-_{\boc_+}$ (resp. $B^-_{\boc_-}$)
consists only of terms in root vectors with positive real part
(resp. negative real part), it follows $\boc_- \neq 0.$ By
\cite[Theorem 5.3 (iii)]{K2} we have $B^-_\boc u_\lam = 0$.  However,
by assumption $B^-_\boc u_\lam \in L(\lam)$ such that $B^-_\boc u_\lam
\neq 0 \mod qL(\lambda)$. This is a contradiction.  
(ii) From Proposition \ref{greatid} we have $\tilde P_{i,-k} u_\lambda
= F_{\delta - \alpha_i}^{(k)} F_{\alpha_i}^{(k)} u_\lambda.$ Since the
weights of $V(\lambda)$ are in the convex hull of $W \lambda$
(\cite[Corollary 5.2]{K2}), this implies that $\tilde P_{i,-k}
u_\lambda = 0$ for $k > \lambda_i$.  Note that for any $i$, $\ell(\rho^{(i)})
\le \lambda_i \iff {\rho'}^{(i)}_1 \le \lambda_i.$ Since the $\tilde
P_{i,-k}$ all commute, considering the top row of the determinant
$S^-_{\boc_0}$, we have $S^-_{\boc_0} u_\lambda = 0$ for $\boc_0
\notin \boc_0(\lambda).$
\end{proof} 

Let $z_i$ be the $\U'$ automorphism of $V(\varpi_i)$ defined in
\cite[\S 5.2]{K2}.
\begin{lem} Let \label{piaction} $i = 1, \dots, n.$   Then on $V(\varpi_i)$:
\begin{equation*}
\tilde P_{i,-1} u_{\varpi_i} = \overline{\tilde P_{i,-1} u_{\varpi_i}}
= z_i^{-1} u_{\varpi_i}, \ \tilde P_{i,-k} u_{\varpi_i} = 0, \ k > 1.
\end{equation*}
\end{lem}

\begin{proof}  For $k > 1$ the statement follows  from 
Proposition \ref{greatid}. Let $k=1$. $V(\varpi_i)$ has a unique
global basis element of weight $\varpi_i - \delta$, which by
definition equals $S_w u_{\varpi_i} = z_i^{-1} u_{\varpi_i},$ where $w
= t(\alpha_i).$ Since $\varpi_i$ is regularly $t(\alpha_i)$--dominant
(see \cite[\S 3.1]{K2}), it follows from the identity (\cite[Appendix
B]{K1}):
$$S_j (b_1 \otimes t_\mu \otimes u_{-\infty}) = {\tilde
f_j}^{a} b_1 \otimes t_\mu \otimes u_{-\infty} \text{ if }
a = \lan h_j, wt(b_1) + \mu \ran \ge 0, $$ that $z_i^{-1} u_{\varpi_i}
\mod q L(\varpi_i) \in B(\infty) \otimes t_{\varpi_i} \otimes u_{-\infty}.$

Since $B^-_\boc u_{\varpi_i} = 0$ for $\boc_- \neq 0,$ the unique
element of the crystal of $B(\varpi_i)$ of weight $\varpi_i - \delta$
must be ${\tilde P}_{i,-1} u_{\varpi_i} \mod q L(\varpi_i).$ Note that
the globalization of an element in $B(\infty) \otimes t_{\varpi_i}
\otimes u_{-\infty}$ remains in $U_q^- u_{\varpi_i}$.  Since for $1
\le j \neq i \le n, \, {\tilde P}_{j,-1}u_{\varpi_i} = 0$, we have
immediately that $G(\tilde P_{i,-1} u_{\varpi_i} \mod qL(\varpi_i)) =
\tilde P_{i,-1} u_{\varpi_i}$, which completes the proof.
\end{proof}

\subsection{The map $\Phi_\lambda$}

Let $\lam = \sum_{i\in I} m_i\varpi_i \in P^0_+$.  The module $V' =
\otimes_{i\in I}V(\varpi_i)^{\otimes m_i}$ has a crystal base
$(L(V'),B(V')) = (L(\varpi_i)^{\otimes m_i}, B(\varpi_i)^{\otimes
m_i}).$ Let $u' = \mathop\otimes\limits_{i\in I}u_{\varpi_i}^{\otimes
m_i}.$
 
For each $i$, and each of the $\nu = 1, \dots, m_i$ factors of
$V(\varpi_i)^{\otimes m_i}$, we let $z_{i,\nu}$ be the commuting
automorphisms defined \cite[\S 4.2]{K2}. By \cite[Theorem 8.5]{K2}, the
submodule
\begin{equation} \label{defW}
W' = U_q(\ag)_\Q[z_{i,\nu}^{\pm 1}]_{1 \le i \le n, 1 \le \nu \le m_i} u'
\subset V'
\end{equation}
has a global crystal base $(L(W'), B(W'))$ such that $L(W') \subset
\otimes_i L(\varpi_i)^{\otimes m_i}$, $B(W') = \otimes_i
B(\varpi_i)^{\otimes m_i}$.  Since $W'$ contains the extremal vector
$u'$ of weight $\lam$ we have a unique $\U$-linear morphism:
\begin{equation} \label{mainmap} \Phi_\lam\colon V(\lam)\to W', 
\end{equation} 
sending $u_\lam$ to $u',$ and which commutes with the crystal
operators $\tilde e_i, \tilde f_i$.

For each $n$--tuple of partition $\boc_0 = (\rho^{(1)},
\rho^{(2)}, \dots, \rho^{(n)})$ we consider the product of Schur
functions in the variables $z^{-1}_{i,\nu}$ (see \cite{M}):
\begin{equation}
s_{\boc_0}(z_{i,\nu}^{-1}) = \prod_{i = 1}^n
s_{\rho^{(i)}}(z_{i,1}^{-1}, \dots, z_{i,m_i}^{-1}).
\end{equation}
Note that for each $i$, $s_{\rho^{(i)}}(z^{-1}_{i,\nu})$ acts as
the $0$ map if $m_i < \ell(\rho^{(i)}).$ We will omit the indices
$i,\nu$ and write $s_{\boc_0}(z^{-1}).$

Using Lemma \ref{piaction} we have:
\begin{prop} \label{schurcor} 
Let $\boc_0 = (\rho^{(1)},\rho^{(2)}, \dots, \rho^{(n)})$ be an
$n$--tuple of partitions:
\begin{equation} \label{schur}
            \Phi_\lam(S^{-}_{\boc_0}\, u_\lambda) =
            s_{\boc_0}(z^{-1})\, u'.
\end{equation}
\end{prop}

\begin{proof} 
Note that $\sigma \circ (\psi \times \psi) \circ \Delta(a) =
\Delta(\psi(a))$ for $a \in U_q$.  Since our $\tilde{P}_{i,-k}$ are
those given in \cite{BCP} after applying $- \circ \psi$ we have by
\cite[Proposition 3.4]{BCP} and \cite{D}
\begin{equation*}
\Delta (\ov{\tilde{P}}_{i,-k}) =\sum_{s=0}^{k} \ov{\tilde{P}}_{i,-s} \ot
\ov{\tilde{P}}_{i,s-k} + {\text{terms acting as }} 0 \text{ on }
v_{\varpi_{j_1}} \otimes v_{\varpi_{j_2}} \text{ for all } j_1, j_2
\in I.  
\end{equation*}
This implies that $\Delta^{m_i}(\ov{\tilde P}_{i,-k})$ acts as
$e_k(z_{i,1}, \dots, z_{i,m_i})$ on $V'$ where $e_k$ is the $k$--th
elementary symmetric function.  Since polynomials in the $\ov {\tilde
P}_{i,-k}$ (resp. elementary symmetric functions) generate the Schur
functions $\ov{S^{-}_{\boc_0}}$ (resp.  $s_{\boc_0}(z^{-1})$) we have
$\Phi_\lam(\ov{S^{-}_{\boc_0}}\, u_\lambda) = s_{\boc_0}(z^{-1})\,
u'.$ Since $\Phi_\lam$ is uniquely defined as a $U_q$ homomorphism, it
commutes with the respective $-$ actions, where $-$ on $L(W')$ is 
$c^{\text{norm}}$ as defined in \cite[\S 8]{K2}.  Now since the
$z_{i,\nu}$ commute with the bar action on $W'$ the proposition
follows.
\end{proof}

Next we consider the image of $B(\lambda)$ under $\Phi_\lambda.$ By
\cite[Theorem 5.1]{K2} every element of $B(\lam)$ is connected to an
extremal vector of the form $b_1 \otimes t_\lam \otimes u_{-\infty},$
which by Proposition \ref{iextremal} equals $S^-_{\boc_0}u_\infty
\otimes t_\lambda \otimes u_{-\infty} \mod q L(\lam).$ Therefore we
have,
\begin{equation} \label{mvtoimaginary}
B(\lambda) = \{X_1\, X_2\, \dots\, X_n (S^-_{\boc_0} u_\lambda\mod q
L(\lambda))\, |\, X_i = \tilde e_i \text{ or } \tilde f_i, \boc_0 \in
\boc_0(\lambda) \} \setminus \{0\}.
\end{equation}
Since $\Phi_{\lambda}$ commutes with crystal operators, and the
$z_{i,\nu}$ induce automorphisms of the $U'_q$--crystal of $V(\varpi_i)$, we
have that $\Phi_\lambda(L(\lambda)) \subset L(W').$ Denote by
$\Phi_{\lambda | q=0}$ the induced map $L(\lam)/qL(\lam) \rightarrow
L(W')/qL(W')$.

\begin{prop} \label{crystalcorresp} Let $B_0(W')$ be the connected
component of $B(W')$ containing $u'.$ Then
$$\Phi_{\lambda | q=0}: \{ b \, | \, b \in B(\lambda) \} \rightarrow
\{s_{\boc_0}(z^{-1}) b'\, |\, \boc_0 \in \boc_0(\lam),\, b' \in B_0(W')
\}$$ is a bijection.
\end{prop}
\begin{proof}
We have
$$\Phi_{\lambda | q=0}(B(\lambda)) \setminus \{0\} = \{
s_{\boc_0}(z^{-1}) B_0(W')\}.$$ Arguing using \eqref{mvtoimaginary},
we check that $\Phi_{\lambda | q=0}(B(\lambda))$ is injective.  Let $b
\in B(\lambda)$ such that $\Phi_{\lambda | q=0} (b) = 0$. Since $b$ is
connected by crystal operators to $b_1 \otimes t_\lambda
\otimes u_{-\infty},$ where $b_1 = S^-_{\boc_0} u_\infty \mod
qL(\infty)$, $\boc_0 \in \boc_0(\lambda)$, this implies
$\Phi_{\lambda | q=0} (S^-_{\boc_0} u_\lambda \mod qL(\lambda)) =
0$. This contradicts Proposition \ref{schurcor}.
\end{proof}

\begin{cor} The map $\Phi_\lambda$ is injective. 
\end{cor}
\begin{proof} 
Since $\Phi_{\lam | q=0}: L(\lam)/qL(\lam) \rightarrow L(W')/qL(W')$
maps the crystal base $B(\lam)$ bijectively, it follows
$\{s_{\boc_0}(z^{-1}) B_0(W')\}$ is linearly independent in
$L(W')/qL(W').$ Write an element $v \in \ker \Phi_\lambda, v \neq 0,$
in terms of the global base $\{G(b)\, | \, b \in B(\lambda)\}$ as $v =
\sum_b c_b(q) G(b).$ Multiplying by a power of $q$ we may assume that
each $c_b(q)$ is regular at $q=0,$ so that $v \mod q L(\lambda) \neq
0$. This implies $\Phi_{\lam | q=0} (v \mod q L(\lam)) \neq 0$, which is
a contradiction.
\end{proof}

By Proposition \ref{crystalcorresp}, for each $b \in B(\lam)$ there
exist $b' \in B_0(W')$ and $s_{\boc_0}(z^{-1})$ such that $\Phi_\lam
(b) = s_{\boc_0}(z^{-1}) b' \mod q L(W')$.  Let $G(b)$, $G(b')$ be the
respective globalizations of $b$ and $b'$.  Then $\Phi_\lam (G(b)) =
s_{\boc_0}(z^{-1}) G(b') \mod q L(W')$.  Since $\Phi_\lam$ commutes
with the $-$ involution, $\Phi_\lam (G(b)) = s_{\boc_0}(z^{-1}) G(b')
\mod q^{-1} \ov{L(W')}$. We conclude:
\begin{thm} $\Phi_\lam$ induces a bijection between the sets
\begin{equation*}
\Phi_\lam : \{G(b) | b \in B(\lam)\} \rightarrow
\{s_{\boc_0}(z^{-1}) G(b)\, | \boc_0 \in \boc_0(\lam), \, b \in B_0(W') \}.
\end{equation*}
\end{thm}

\begin{rem} Taken together the results of this section give the
conjectures \cite[13.1, 13.2]{K2}. To obtain 13.1 (iii), consider that
the crystal $\bigotimes_{i\in I} B(m_i\varpi_i)$ is by Proposition
\ref{crystalcorresp} in bijective correspondence with
$\{s_{\boc_0}(z^{-1}) B_0(W')\}$, and note that $\Phi_\lambda$ factors
through $\bigotimes_{i\in I} V(m_i\varpi_i)$.
\end{rem}

\end{document}